\documentclass[a4paper]{article}
\usepackage[english]{babel}
\usepackage{amsthm}  % for math
\usepackage{amsmath} % for math
\usepackage{amssymb} % for more math
\usepackage{mathabx} % for even more math
\usepackage{graphicx}% for graphicx and figures
\usepackage{fancyhdr}% for the footer
\usepackage{setspace}% for stuff like doublespacing
\usepackage{lineno}	 % for the linenumbers
\usepackage{lastpage}% for finding the last page (needed in the footer)
\usepackage{color}
\usepackage{tikz}
\usepackage[percent]{overpic}
\usepackage{multirow}
\usepackage{array}
\usepackage{bigdelim}
\usepackage{subcaption}

\definecolor{yelloworange}{RGB}{255,204,0}
\definecolor{forestgreen}{rgb}{0.13,0.55,0.13 }
\definecolor{magenta}{RGB}{255,0,255}
\definecolor{cerulean}{rgb}{0,0.48,0.65}

\newcommand{\mytitle}{Optimal Deterministic Algorithm Generation }
\newcommand{\myshorttitle}{Optimal Deterministic Algorithm Generation}
\newcommand{\myauthor}{\textbf{Alexander Mitsos, Jaromi{\l} Najman} \\
  AVT Process Systems Engineering (SVT), RWTH Aachen University, Turmstrasse 46, Aachen, 52064, Germany, amitsos@alum.mit.edu ORCID: 0000-0003-0335-6566\\
 \textbf{Ioannis G. Kevrekidis} (corresponding, yannis@princeton.edu)  \\
Department of Chemical and Biological Engineering \& Program in Applied and Computational Mathematics, Princeton University, A319 Engineering Quad;
Princeton, NJ 08544, yannis@princeton.edu, Phone: 609-258-2818 (also IAS-TUM, Garching and Zuse Institut, FU Berlin, Germany.)} % myauthor is important for the ^b !!!
\newcommand{\myauthorshort}{A. Mitsos et al.}
\author{\myauthor}
\usepackage[colorlinks,linkcolor=blue,citecolor=blue,urlcolor=blue,pdftitle={\mytitle},pdfauthor={A. Mitsos, J. Najman, I.G. Kevrekidis}]{hyperref} 

{
	\fancyhead[L]{\footnotesize{\myshorttitle}}
	\fancyhead[R]{\footnotesize{\textit{\the\day.\the\month.\the\year}}}
	\fancyfoot[L]{\copyright \footnotesize{\textit{\myauthorshort}}}
	\fancyfoot[R]{\footnotesize{Page \thepage\ of \pageref*{LastPage}}}
	\cfoot{}
}

\fancypagestyle{firststyle}
{
	 \fancyhead[L]{\copyright \footnotesize{\textit{\myauthorshort}}}
	 \fancyhead[C]{}
   \fancyhead[R]{\footnotesize{Page \thepage\ of \pageref*{LastPage}}}
	 \fancyfoot[L]{}
   \fancyfoot[C]{}
	 \fancyfoot[R]{}
}

\addto\captionsenglish{%
}

\newtheorem{Definition}{Definition}
\theoremstyle{remark} % to distinguish between different theorem styles (this is without italic)

\makeatletter
\let\@addpunct\@gobble
\g@addto@macro{\thm@space@setup}{\thm@headpunct{}} % point behind Def., Proof etc.
\makeatother

\renewenvironment{abstract}{\noindent\textbf{Abstract:}}{}
\newenvironment{keywords}{\noindent\textbf{Keywords:}}{}
\newenvironment{AMS}{\noindent\textbf{AMS:}}{}

\begin{document}

\thispagestyle{firststyle}
\begin{flushleft}\begin{large}\textbf{\mytitle}\end{large} \end{flushleft}
\begin{small}\myauthor\end{small} 

\vspace{0.5cm}

\begin{abstract}
 A formulation for the automated generation of algorithms via mathematical programming (optimization) is proposed.
  The formulation is based on the concept of optimizing within a parameterized family of algorithms, or equivalently a family of functions describing the algorithmic steps.
  The optimization variables are the parameters  -within this family of algorithms-  that encode algorithm design: the computational steps of which the selected algorithms consists.
  The objective function of the optimization problem encodes the merit function of the algorithm, e.g., the computational cost (possibly also including a cost component for memory requirements) of the algorithm execution.
  The constraints of the optimization problem ensure convergence of the algorithm, i.e., solution of the problem at hand.
  The formulation is described prototypically for algorithms used in solving nonlinear equations and in performing unconstrained optimization; the parametrized algorithm family considered is that of monomials in function and derivative evaluation (including negative powers).
  A prototype implementation in GAMS is developed along with illustrative results demonstrating cases for which well-known algorithms are shown to be optimal.
  The formulation is a mixed-integer nonlinear program (MINLP).
  To overcome the multimodality arising from nonconvexity in the optimization problem, a combination of brute force and general-purpose deterministic global algorithms is employed to guarantee the optimality of the algorithm devised.
  We then discuss several directions towards which this methodology can be extended, their scope and limitations.
\end{abstract}  

\begin{keywords} Optimization, Nonlinear Equations, Algorithms, Optimal Control \end{keywords}

\begin{AMS}
49M05, 49M15, 49M25, 49M37, 49L20, 90C90, 90C26
\end{AMS}

\section{Introduction}
Computation has led to major advances in science and engineering, in large part due to ingenious numerical algorithms.
The development of algorithms is thus of crucial importance and requires substantial resources and time.
Typically multiple algorithms exist for a given task.
Comparisons of algorithms that tackle a given family of problems is sometimes performed theoretically and sometimes numerically.
It would be desirable to automatically generate algorithms and have a guarantee that these are the best for a given problem or for a class of problems.
We propose the use of numerical optimization to design \emph{optimal algorithms}, i.e., algorithms that perform better than any conceivable alternative.
More precisely, each iteration of the algorithms is interpreted, in an input-output sense, as a function evaluation. Then, an optimization problem is formulated that finds the best (in a given metric) algorithm among a family of algorithms/functions.

There is large literature on automatic generation of algorithms.
In many cases genetic algorithms are used to first generate a set of algorithms built of elementary components of promising approaches. The set is evaluated on test problems and a new elitist algorithm is obtained by combination \cite{bain2004methods,koza1992genetic,koza1994genetic}. The work by \cite{koza1992genetic,koza1994genetic} covers basic methods and ideas in the field of genetic programming. Automatic algorithm generation for retrieval, classification and other data manipulation has been proposed by \cite{kuhner2002automatic} based on statistical data analysis and other methods.
For instance, in \cite{parada2013automatic} algorithms for the well-known knapsack problem are automatically generated by genetic algorithms.
Algorithms are first generated for small knapsack instances and then compared to algorithms computed on larger test sets. %In \cite{parada2013automatic} an increase in efficiency is noted if larger instantiations for the algorithm generation are used.
In \cite{garber2008automatic} work is presented on automatic generation of parallel sorting algorithms. They combine known sorting algorithms to obtain a best performing algorithm for a certain input.
The well known art gallery problem (AGP) is discussed by \cite{tozoni2013practical}, where an iterative algorithm is developed and in that sense resembles our proposal. %The problem considered is to find optimal positions for guards to cover the art gallery completely, which is a given polygon. Although the problem at first sounds simple, its combinatorial nature explosively increases its complexity. The algorithm presented by \cite{tozoni2013practical} is iterative, which in some sense resembles the proposed idea of optimal algorithm generation through optimization.
Performance analysis of first-order methods for smooth unconstrained convex optimization can be of high interest regarding our work if new first-order algorithms are found via our approach. Building on the work of \cite{nemirovsky1983problem}, in \cite{Drori2014} a novel approach for performance analysis is presented. In \cite{Kim2016} Nesterov-like first-order algorithms for smooth unconstrained convex optimization are developed. The worst-case convergence bound of the presented algorithms is twice as small as that of Nesterov's methods.  
Moreover, there is substantial work on automatic code and  algorithm generation for certain problems, e.g., \cite{ricart1981optimal,arya1994optimal,bacher1996automatic,coelho1999robust}.
While this approach is fundamentally different, in a sense it resembles the proposal herein: both require an external algorithm in order to compute a sub-result. Many depend on algebraic operations regarding matrices, e.g., QR decomposition, diagonal transformations or eigenvalue computation. Many of those operations are well documented in the Thesis of \cite{bientinesi2006}.
Finally, the problem of tuning the parameters of existing algorithms through optimization has been considered~\cite{weinan2016}.
Lessard et al.~\cite{lessard_16_1} use control theory to derive bounds on the convergence rate of important algorithms, most of which are considered herein as well.

Consider an iterative algorithm, like the time-honored Newton-Raphson for finding roots of algebraic equations, e.g., in one dimension
$ x_{n+1} = x_n - f(x_n)/f'(x_n) \equiv x_n + u_n(x_n). $
The iterates can thus be seen as a ``control action'' $u_n$ that depends only on the current guess, i.e., ``state-dependent feedback''; a similar interpretation can be given to most algorithms.
%In the same spirit, we can probably write most usual algorithms as feedback laws - the update is the ``state-dependent feedback'' on the state itself.
%One can have variants (like Nesterov's method) where the feedback depends on the last two iterations, or where the feedback is continuous in time (consider the two recent beautiful embeddings of Newton as two-dimensional ordinary differential equations by Candes and coworkers and by Jordan and coworkers).
But if algorithms are feedback laws (towards a prescribed objective), then we do have celebrated mathematical tools for computing optimal feedback, e.g., Hamilton-Jacobi-Bellman. It would thus make sense to use these tools to devise optimal algorithms.

The key goal of the manuscript is to formulate a systematic process for obtaining optimal algorithms.
To the best of our knowledge, there exists no other proposal in the literature to use deterministic global optimization for finding an algorithm that is optimal w.r.t. a certain objective. 
Optimality is determined based on a desired measurable performance property which is encoded as the objective of the optimization problem. For instance this can be the computational expense of the algorithm, possibly with the memory requirements weighting in.
Upon convergence of our formulation, the resulting algorithm is optimal among a considered family of algorithms which in turn is encoded using the variables (or degrees of freedom) of the optimization problem.
The optimization formulation is completed by the constraints which ensure that only feasible algorithms will be selected among the family of algorithms considered, i.e., algorithms that converge to a solution of the problem at hand.
The formulations can be cast as nonlinear programs (NLP) or mixed-integer nonlinear programs (MINLP).
These formulations are quite expensive; in fact, as we discuss in the following they are expected to be more expensive than the original problem; however, we argue that this is acceptable.

Our proof of concept illustrations involves devising optimal algorithms that perform local unconstrained minimization of scalar functions and those that solve equation systems.
The class of algorithms within which optimality is sought herein involve evaluations of the function and its first two derivatives; they can be composed of up to ten iterations of a monomial involving these quantities, what we call a monomial-type algorithm.
The class is then extended to allow for methods such as Nesterov's acceleration.
The well-known steepest descent is shown to be optimal in finding a stationary point of a univariate fourth-order polynomial, whereas it is infeasible for the two-dimensional Rosenbrock objective function (it would take more than the allowed number of iterations to converge). Both statements hold for a given initial guess (starting value of the variables given to the algorithm).
We also consider the effect of different initial guesses and show that Nesterov's method is the cheapest in an ensemble of initial guesses.
Finally, we focus on finding optimal algorithms for a given problem, e.g., minimizing a fixed function; however, we discuss how to handle multiple instances by considering families of problems.

The remainder of the manuscript is structured starting with the simplest possible setting, namely solution of a single nonlinear equation by a particular class of algorithms. Two immediate extensions are considered (a) to  equations systems and (b) to local unconstrained optimization.
Subsequently, numerical results are presented along with some discussion.
Potential extension to more general problems and algorithm classes are then discussed as well as limitations and relation to computational complexity.
Finally, key conclusions are presented.

\section{Definitions and Assumptions}
For the sake of simplicity, first, solution of a single nonlinear equation is considered as the prototypical task an algorithm has to tackle.
\begin{Definition}[Solution of equation]
  \label{def:eqnsol}
  Let $x \in X \subset \mathbb{R}$ and $f:X \rightarrow \mathbb{R}$.
  A solution of the equation is a point $x^* \in X$ with $f(x^*)=0$.
  An approximate solution of the equation is a point $x^* \in X$ with $|f(x^*)|\leq \varepsilon$, $ \varepsilon > 0$.
\end{Definition}
Note that no assumptions on existence or uniqueness are made, nor convexity of the function.
A particular interest is to find a steady-state solution of a dynamic system $\dot{x}(t)=f\left(x(t)\right)$.

We will first make use of relatively simple algorithms:
\begin{Definition}[Simple Algorithm]
  Algorithms operate on the variable space ${\bf x} \in X \subset \mathbb{R}^{n_x}$ and start with an initial guess ${\bf x}^{(0)} \in X$.
  A given iteration $it$ of a simple algorithm amounts to calculating the current iterate as a function of the previous ${\bf x}^{(it)}=g_{it}\left({\bf x}^{(it-1)}\right)$, $g_{it}: X \rightarrow X$.
  Therefore, after iteration $it$, the algorithm has calculated ${\bf x}^{(it)}=g_{it}\left(g_{it-1}\left(\dots(g_{2}(g_{1}({\bf x}^{(0)})))\right)\right)$.
  A special case are algorithms that satisfy $g_{it_1}({\bf x})=g_{it_2}({\bf x})$ for all $it_1,it_2$, i.e., use the same function at each iteration.
\end{Definition}
Note that the definition includes the initial guess ${\bf x}^{(0)}$. Herein both a fixed initial guess and an ensemble of initial conditions is used.
Unless otherwise noted the same function will be used at each iteration $g_{it}=g$.

In other words, algorithmic iterations can be seen as functions and  algorithms as composite functions.
This motivates our premise, that finding an optimal algorithm constitutes an optimal control problem.
In some sense the formulation finds the optimal function among a set of algorithms considered.
One could thus talk of a ``hyper-algorithm'' or ``meta-algorithm'' implying a (parametrized) family of algorithms, or possibly the span of a ``basis set'' of algorithms that are used to identify an optimal algorithm.
%\begin{Assumption} %[Host set]
In the following we will assume $X=\mathbb{R}^{n_x}$ and consider approximate feasible algorithms:
%\end{Assumption}
%The goal of the proposed approach is to automatically generate algorithms that solve a problem in some optimal fashion.
%This motivates the following definition.
\begin{Definition} %[Feasible and Optimal Algorithms]
  An algorithm is \emph{feasible} if it solves a problem; otherwise it is termed \emph{infeasible}.
  More precisely: it is termed \emph{feasible in the limit} for a given problem if it solves this problem as the number of iterations approaches $\infty$; it is termed \emph{finitely feasible} if it solves the problem after a finite number of iterations; it is termed \emph{approximate feasible} if after a finite number of iterations it solves a problem approximately. In all these cases the statements holds for some (potentially open) set of initial conditions.
  A feasible algorithm is \emph{optimal} with respect to a metric if among all feasible algorithms it minimizes that metric.
\end{Definition}
One could also distinguish between feasible path algorithms, i.e., those that satisfy ${\bf x}^{(it)} \in X$, $\forall it$ in contrast to algorithms that for intermediate iterations violate this condition.
%Herein, unless otherwise noted, we will consider approximate feasible algorithms.

%\begin{Definition} %[Notation of Derivative]
%  Consider $f:X \subset \mathbb{R} \rightarrow \mathbb{R}$.
%  The derivative of order $j$ is denoted by $f^{(j)}$ with $f^{(0)}\equiv f$.
%\end{Definition}

%\begin{Assumption} %[Smoothness]
Throughout the report it is assumed that $f:X \subset \mathbb{R} \rightarrow \mathbb{R}$ is sufficiently smooth in the sense that if algorithms are considered that use derivatives of a given order, then $f$ is continuously differentiable to at least that order. The derivative of order $j$ is denoted by $f^{(j)}$ with $f^{(0)}\equiv f$.
%\end{Assumption}
A key point in our approach is the selection of a class of algorithms: a (parametrizable) family of functions.
It is possible, at least in principle, to directly consider the functions $g_{it}$ and optimize in the corresponding space of functions.
Alternatively, one can identify a specific family of algorithms, e.g., those based on gradients, which includes well-known algorithms such as Newton and explicit Euler:
\begin{Definition} %[Monomial-Type Algorithms]
\label{def:ratpolalgo}
  Consider problems that involve a single variable $x \in X \subset \mathbb{R}$ and $f: X \rightarrow \mathbb{R}$.
  \emph{Monomial-type algorithms} are those that consist of a monomial $g_{it}(x)=x+ \alpha_{it} \Pi_{j=0}^{j_{max}} \left(f^{(j)}(x)\right)^{\nu_j}$ of derivatives of at most order $j_{max}$, allowing for positive and negative (integer) powers $\nu_j$.
\end{Definition}
%The class of monomial-type algorithms includes well-known algorithms such as Newton or explicit Euler.
As examples of algorithms in this family consider the solution of a single nonlinear equation, Definition~\ref{def:eqnsol}, by the following important algorithms
\begin{itemize}
\item Newton: $x^{(it)}=x^{(it-1)}+\frac{f^{(0)}(x^{(it-1)})}{f^{(1)}(x^{(it-1)})}$, i.e., we have $\alpha=1$, $\nu_0=1$, $\nu_1=-1$, $\nu_{i>1}=0$.
\item Explicit Euler: $x^{(it)}=x^{(it-1)}+f(x^{(it-1)}) \Delta t$, i.e., we have $\alpha=\Delta t$, $\nu_0=1$, $\nu_{i>0}=0$.
\end{itemize}

\section{Development}
\label{sec:formu}
We first develop an optimization formulation that identifies optimal algorithms of a particular problem: the solution of a given nonlinear scalar equation, Definition~\ref{def:eqnsol}.
We will consider approximate feasible solutions and simple algorithms that use the same function at each iteration.
If the optimization metric accounts for the computational cost of the algorithm, then it is expected that algorithms that are only feasible in the limit will not be optimal since their cost will be higher than those of finitely feasible algorithms.

Obviously to solve the optimization problem with local methods one needs an initial guess for the optimization variables, i.e., an initial algorithm.
Some local optimization methods even require a ``feasible'' initial guess, i.e., an algorithm that solves (non-optimally) the problem at hand.
Note however, that we will consider deterministic global methods for the solution of the optimization problem, and as such, at least in theory, there is no need for an initial guess.

\subsection{Infinite Problem}
\label{sec:formuinf}
We are interested in devising algorithms as the solution of an optimization problem.
The key idea is to include the iterates of the algorithm $x^{(it)}$ as (intermediate) \emph{optimization variables} that are determined by the algorithmic iteration.
By imposing a maximal number of iterations, we have a finite number of variables.
The algorithm itself is an optimization variable, albeit an infinite one (optimal control).
Moreover, an end-point constraint ensures that only feasible algorithms are considered.
Finally, an optimization metric is required which maps from the function space to $\mathbb{R}$, such as the computational cost.
We will make the assumption that the objective is the sum of a cost at each iteration which in turn only depends on the iterate, and implicitly on the problem to be solved.

\begin{align}
&\min_{g, {x}^{(it)},it_{con}} \sum_{it=1}^{it_{con}} \phi\left(g\left(x^{(it)}\right)\right) \notag \\
&\text{s.t. } x^{(it)}=g\left(x^{(it-1)}\right), \quad it=1,2,\dots, it_{con} \notag\\
&    \qquad        \left|f\left(x^{(it_{con})}\right)\right|\leq \epsilon, \label{optformuinf}
\end{align}
where $it_{con}$ corresponds to the iteration at which the desired convergence criterion is met.
The optimal point found embodies an optimal algorithm $g^*$; depending on the method used it may be a global, an approximately global or a local optimum.
For notational simplicity, we have assumed that the minimum is attained.
It is outside the scope of this manuscript to determine if the infimum is attained, e.g., if $g$ can be described by a finite basis of a finite vector space.

The intermediate variables can be considered as part of the optimization variables along with the constraints or eliminated.
Note the analogy to full discretization vs.~late discretization vs.~multiple shooting in dynamic optimization~\cite{biegler_book_2,barton_06_1} with known advantages and disadvantages.

Formulation \eqref{optformuinf} is not a regular optimal control problem as the number of variables is not a priori known. There are however several straightforward approaches to overcome this problem, such as the following.
A maximal number of iterations $it_{max}$ is considered a priori and dummy variables for $it_{con}<it<it_{max}$ are used along with binary variables capturing if the convergence has been met.
The cost is accounted only for the iterates before convergence:
  \begin{align}
    &\min_{g, {x}^{(it)}} \sum_{it=1}^{it_{max}} y_{res}^{(it)} \phi\left(g\left(x^{(it)}\right)\right) \notag \\
    &\text{s.t. } x^{(it)}=g\left(x^{(it-1)}\right), \quad it=1,2,\dots, it_{max} \\
    & \qquad y_{res}^{it}= \begin{cases} 0 & \text{if } \left|f\left(x^{(it)}\right)\right| < \varepsilon  \\ 1 & \text{else} \end{cases} \\
    & \qquad y_{res}^{it_{max}}=0,
    \end{align}
where satisfaction of the last constraint ensures convergence.
In the Appendix an alternative is given.
In the formulation above each iteration uses the same function $g$. Allowing different functional forms at each iteration step, is of interest, e.g., to devise hybrid iteration algorithms, such as few steepest descent steps followed by Newton steps.

\subsection{Finite Problem}
\label{sec:optformufin}
We now use the family of monomial-type algorithms, Definition~\ref{def:ratpolalgo}, as our ensemble of algorithms.
In analogy to the infinite case we assume that the objective depends on the evaluation of derivatives (of the result of the algorithm with respect to its inputs), which in turn only depends on the iterate, and implicitly on the problem to be solved.
To account for the expense of manipulating the derivatives, e.g., forming the inverse of a matrix (recall that simple algorithms scale cubically in the size, i.e., the number of rows of the matrix), we assume that the cost also involves an exponent:
\begin{align}
&\min_{\nu_j,\alpha_{it}, {x}^{(it)},it_{con}} \sum_{it=1}^{it_{con}} \sum_{j=0}^{j_{max}} \phi_j\left(f^{(j)}\left(x^{(it)}\right),\alpha_{it},\nu_j \right) \notag \\
&\text{s.t. } x^{(it)}=x^{(it-1)}+ \alpha_{it} \Pi_{j=0}^{j_{max}} \left(f^{(j)}(x^{(it-1)})\right)^{\nu_j}, it=1,\dots, it_{con} \notag\\
&   \qquad         \left|f\left(x^{(it_{con})}\right)\right|\leq \epsilon, \label{optformufin}
\end{align}
where again the number of variables is not known a priori, see above.
Note that, since the exponents are integer-valued, we have an MINLP and thus, for a fixed $f$ we expect that the minimum is attained.
Again the intermediate variables can be considered as part of the optimization variables along with the constraints, or they can be eliminated.
The optimal point found gives optimal step size $\alpha^*$ as well as exponents $\nu_j^*$.
Herein, we will restrict the step size to discrete values, mimicking line-search methods, but in general step size can be optimized over a continuous range.
Allowing for a different algorithm found at each step ($\nu_{j}^{it_1} \neq \nu_{j}^{it_2}$), the number of combinations dramatically increases but there is no \emph{conceptual} difficulty.
The form is not common for an MINLP but can be easily converted to a standard one as shown in the Appendix.

It may not always be easy to estimate the cost $\phi_j\left(f^{(j)}\left(x^{(it)}\right),\alpha_{it},\nu_j \right)$.
For instance consider Newton's method augmented with a line-search method
$ x^{(it)}=x^{(it-1)}+\alpha_{it} \frac{f(x^{(it-1)})}{f^{(1)}(x^{(it-1)})}$.
For each major iteration, multiple minor iterations are required, with evaluations of $f$.
If the cost of the minor iterations is negligible, the computational cost is simply the evaluation of $f$, its derivative and its inversion.
If we know  the number of minor iterations then we can calculate the required number of evaluations.
If the number of iterations is not known, we can possibly obtain it by the step size, but this requires knowledge of the line search strategy.
In some sense this is an advantage of the proposed approach: the step size is not determined by a line search method but rather by the optimization formulation.
However, a challenge is that we need to use a good cost for the step size, i.e., incorporate in the objective function the computational expense associated with the algorithm selecting a favorable step size $\alpha_{it}$.
Otherwise, the optimizer can select an optimal $\alpha_{it}$ for that instance which is not necessarily sensible when the cost is accounted for.
To mimic the way conventional algorithms work, in the numerical results we allow discrete values of the step size $\alpha_{(it)}=\pm 2^{-\bar{\alpha}_{it}}$ and put a bigger cost for bigger $\bar{\alpha}$ since this corresponds to more evaluations of the line search.

The algorithm is fully determined by selecting $\nu_j$ and $\bar{\alpha}_{it}$.
So in some sense the MINLP is a purely integer problem: the variables $x^{(it)}$ are uniquely determined by the equations.
To give a sense of the problem complexity, consider the setup used in the numerical case studies. %, Section~\ref{sec:numerical}.
Assume we allow derivatives of order $j \in \{0,1,2\}$ and exponents $\nu_j \in \{-2,-1,0,1,2\}$ and these are fixed for each iteration $it$.
This gives $5^3=125$ basic algorithms.
Further we decide for each algorithm $it$ if the step size $\alpha_{it}$ is positive or negative, and allow $\bar{\alpha}_{it} \in \{0,1,\dots,10\}$, different for each iteration $it$ of the algorithm.
Thus we have $2\times 10^5$ combinations of step sizes for each algorithm and $25 \times 10^6$ total number of combinations.

\subsection{Equation Systems}
\label{sec:eqnsystems}
A natural extension of solving a single equation is to solve a system of equations.
\begin{Definition} %[Solution of equation systems]
  \label{def:eqnssol}
  Let ${\bf x} \in X \subset \mathbb{R}^{n_x}$ and ${\bf f}:X \rightarrow \mathbb{R}^{n_x}$.
  A solution of the equation system is a point ${\bf x}^* \in X$ with ${\bf f}({\bf x}^*)={\bf 0}$.
  An approximate solution of the equation is a point ${\bf x}^* \in X$ with $||{\bf f}\left({\bf x}^*\right)||\leq \varepsilon$,  $\varepsilon>0$.
\end{Definition}
In addition to more cumbersome notation, the optimization problems become much more challenging to solve due to increased number of variables and also due to the more expensive operations.
Obviously the monomials need to be defined appropriately; for instance the inverse of the first derivative corresponds to the inverse of the Jacobian, assuming it has a unique solution. This in turn can be written as an implicit function or an equation system can be used in the optimization formulation. For the numerical results herein we use small dimensions (up to two herein) and analytic expressions for the inverse.
In the Appendix we discuss an approach more amenable to higher dimensions.

\subsection{Special Two-Step Methods}
A particularly interesting class of algorithms are the so-called Nesterov's methods, see e.g., \cite{candes_15_1}.
These compute an intermediate iterate $y^{(it)}=h(x^{(it)},x^{(it-1)})$ and use it for the computation of the next value $x^{(it)}$.
The above formulation does not allow such intermediate calculations but can be relatively easily extended to 
\begin{alignat}{2}
&\min_{g, {x}^{(it)},it_{con}} &&\sum_{it=1}^{it_{con}} \phi\left(g\left(x^{(it)}\right)\right) \notag \\
&\quad\quad\text{s.t.} &&x^{(it)}=g\left(y^{(it-1)}\right), \quad it=1,2,\dots, it_{con} \notag\\
&             &&y^{(it)}=x^{(it)}+\beta_{it}\left(x^{(it)}-x^{(it-1)}\right), \quad it=1,2,\dots, it_{con} \notag\\
&            &&\left|f\left(x^{(it_{con})}\right)\right|\leq \epsilon \label{Nesterov}.
\end{alignat}
This latter formulation covers most of Nesterov's constant schemes along with potentially new schemes.
The introduction of the parameter $\beta$ necessitates capturing the cost of different values of $\beta$.
Herein, we add the value of $\beta$ to the term of each iteration for the objective function.
We note that this choice is somewhat questionable but for the illustrative, proof-of-concept computations herein, is acceptable.

\subsection{Finding Optimal Local Optimization Algorithms}
It is relatively easy to adapt the formulation from equation solving to local solution of optimization problems.
The only substantial change consists in replacing the end-point constraint with some desired termination criterion such as an approximate solution of the KKT conditions.
For unconstrained problems the corresponding criterion is stationarity of the objective $\frac{\partial f}{\partial x_{ix}}=0$.
Similarly to equation systems increased dimensionality ($n_x>1$) implies that vectors and matrices arise and operations need to be appropriately formulated.

\section{Numerical Results}

\label{sec:numerical}

\subsection{Method}
We develop a prototype implementation in GAMS~\cite{gams_ref}, inspired by \cite{mitsos_08_siprestricted}. % using also the ideas of Michael Bussieck for the use of macros.
We used GAMS 24.5 and tested both local and global methods, mostly KNITRO~\cite{knitro_ref} and BARON~\cite{tawarmalani_05_1}.
%All the common elements are placed in ``optforalgomain.gms''.
%The problem-specific definitions are then given in corresponding GAMS files that are included so it is easy to switch between problems.
To limit the number of variables we avoid introducing variables for $f$ and its derivatives and rather define them using macros.
We do however introduce auxiliary variables for the iterates of the algorithm as well as variables that capture the residuals, see below.
Both could be avoided by the use of further macros.
Since we utilize global solvers we impose finite bounds for all variables.

As discussed, we encode the choice of derivative as a set of binary variables $y_j^k$.
\[  \sum_{k=k_{min}}^{k_{max}} y_j^k \left(f^{(j)}(x^{(i-1)})\right)^k \qquad \sum_{k=k_{min}}^{k_{max}} y_j^k=1.  \]
Similarly we encode the step size using discrete values $\pm 1/2^{\alpha_{it}}$.
Two problems are considered, namely solution of equation systems and minimization, and in both cases we consider one-dimensional and two-dimensional problems.
The minimization is assumed unconstrained, but since we impose bounds on the variables we essentially assume knowledge of an open box in which optimal solutions are found.
%For two-dimensional problems, for simplicity, we do not allow second derivatives.
We allow 10 steps of the algorithm. We require that at least the last iterate meets the desired convergence criteria.

For each iterate $it$ we calculate the residuals by introducing a real-valued variable $res^{(it)}$ and a binary $y_{res}^{(it)}$ denoting if the desired convergence is met ($y_{res}^{(it)}=0$) or not ($y_{res}^{(it)}=1$).
For the case of equation solving we take as residuals the value of the right-hand-side of the equation (function) $res^{it}=\max_{if} |f_{if}\left({\bf x}^{(it)}\right)|$ and in the case of minimization the value of the derivative $res^{it}=\max_{ix} \left|\left.\left(\frac{\partial f}{\partial x_{ix}}\right)\right|_{{\bf x}^{(it)}}\right|$.
This corresponds to the infinity norm and this is chosen since initial testing suggests more robust optimization.
The absolute values are written as pair of inequalities, e.g., for equation solving
\[ res^{(it)} \geq \pm f_{if}\left({\bf x}^{(it)}\right) \forall if. \]
  Then the binary variable is decided based on a standard big-M constraint
\[ M_{res} y_{res}^{(it)} \geq (res^{(it)} -\varepsilon), \]
where $M_{res}$ is the big-M constant that captures the magnitude of the residual $res$.
Only iterations that do not meet the residual tolerance are counted for in the objective
\[ \min \sum_{it=1}^{it_{max}} y_{res}^{(it)} \left( \bar{\alpha}_{it} +\sum_{j=1}^{j_{max}} \sum_{k=k_{min}}^{k_{max}} y_j^k  \phi_j^k \right).  \]
For the cost constants $\phi_j^k$ we assume a zero cost for zero exponent, a cost of 1 for exponent 1, a cost of 1.5 for an exponent of 2, a cost of 2 for an exponent of -1 and a cost of 3 for an exponent of -2.
The idea is to capture the expense of inversion and exponentiation.
These costs are multiplied by 1 for $f$, 10 for $f^{(1)}$ and 100 for $f^{(2)}$.

The problems considered are
\begin{itemize}
  \item Solving $x^3=1$, $x \in [-2,2]$, $x^{(0)}=0.1$
  \item Solving $x \exp(x)=1$, $x \in [-2,2]$, $x^{(0)}=0.1$
  \item $\min_x x^4+x^3-x^2-1$,  $x \in [-2,2]$, $x^{(0)}=0.1$
  \item Rosenbrock function in 2d $\min_{{\bf x}} 100 (x_2-x_1^2)^2+(1-x_1)^2$,   ${\bf x} \in [-2,2]^2$, ${\bf x}^{(0)}=(0.7,0.75)$
  \item Simple 2d equation solving $x_2-x_1^2=0$,  $5x_2-\exp(x_1)=0$. Here we  excluded second derivatives.
  \item $\min_x \exp(x)+x^2$ (strongly convex with convexity parameter $2$), $x \in [-2,2]$, $x^{(0)}=0.1$
  \item $\min_{\bf x} (x_1-1)^2+2x_2^2-x_1x_2$, ${\bf x} \in [-1,2]^2$, ${\bf x}^{(0)}=(-1,-1)$
  \item $\min_{\bf x} \exp(x_1)+x_1^2+\exp(x_2)+x_1^2+x_1 x_2$, ${\bf x} \in [-1,1]^2$, $\mathbf{x}^{(0)}=(0.5,0.5)$. The optimimum of $f$ is at $\mathbf{x}^*\approx (-0.257627653,-0.257627653)$    
\end{itemize}

\subsection{Single-Step Methods}
We tried several runs, mostly with short times (order of minute) but some with longer.
The key findings is that KNITRO often fails to find feasible algorithms, even when they exist.
BARON often also fails to find feasible algorithms within few minutes.
In several cases BARON manages to find feasible algorithms and converge the lower bound but in many cases the lower bound stays at low levels, i.e., the certificate of optimality is relatively weak.
Often the convergence behavior of BARON is improved if it is initialized with KNITRO.
In some cases convergence of the optimization was improved when we enforced that residuals are decreasing with iteration $res^{it}<res^{it-1}$.
It is noteworthy that in early computational experiments the optimizer also managed to furnish unexpected solutions that resulted in refining the formulation, e.g., choosing a step size $\alpha$ that resulted in direct convergence before discrete values of $\alpha$ were enforced.

Due to the relatively difficult convergence of the optimization and the relatively small number of distinct algorithms (125), we solved each problem for each algorithm in three sequential steps:
first KNITRO to find good initial guesses, then BARON without convergence criterion (i.e., without imposing $y_{res}^{it_{con}}=0$) and then BARON with convergence criterion.
Each of the runs was given 60 seconds in a first computational experiment and 600 seconds in a second.
With the small CPU time some runs converged to the specified optimality tolerance of 0.1 (both absolute and relative), while most did not;
with the higher CPU more algorithms were proven to be infeasible while, for some, feasible solutions were found.
Convergence to the global optimum is not guaranteed for all cases.
The main findings of this explicit enumeration is that most algorithms are infeasible.
For equation-solving, in addition to Newton's algorithm some algorithms with second derivatives are feasible.
During the loop some runs resulted in internal errors. %, but since we run via wine we are unsure where the mistakes happened.

The algorithms discovered as optimal for unconstrained minimization are instructive. In the problem $\min_{x \in [-2,2]} x^4+x^3-x^2-1$,  with $x^{(0)}=0.1$ the steepest descent algorithm is optimal.
It is interesting to note that the algorithms do not furnish a global minimum: in the specific formulation used optimal algorithms are computationally cheap algorithms that give a stationary point.
In contrast, for the minimization of the Rosenbrock function in 2d $\min_{{\bf x} \in [-2,2]^2} 100 (x_2-x_1^2)^2+(1-x_1)^2$,   with $x_1^{(0)}=0.7$, $x_2^{(0)}=0.75$ steepest descent is not feasible within ten iterations which is well-known behavior.
In the case of solving a single equation several algorithms were furnished as feasible.
%but these rely on choosing a spurious sequence of steps; in other words these are not interesting algorithms but rather artifacts of the formulation.

\subsubsection{Unexpected Algorithms Discovered as Optimal}
Two particularly interesting algorithms identified (with $\alpha_{it} <0$) are
\begin{equation}
x^{(it)}=x^{(it-1)}+\alpha_{it}\frac{f^{(0)}(x^{(it-1)})f^{(1)}(x^{(it-1)})}{f^{(2)}(x^{(it-1)})}, \label{algo1}
\end{equation}
and
\begin{equation}
x^{(it)}=x^{(it-1)}+\alpha_{it}\frac{f^{(0)}(x^{(it-1)})(f^{(1)}(x^{(it-1)}))^2}{(f^{(2)}(x^{(it-1)}))^2} \label{algo2}
\end{equation}
which can be seen as a combination of Newton's algorithms for equation solving and unconstrained optimization.
The step sizes found vary from iteration to iteration, but the algorithms converge also with $\alpha=-1$.  % It would be interesting to confirm optimality through a HJB approach.

There is an algorithm mentioned between the lines of section 3.1.2 in \cite{deuflhard2004newton} for root finding of \emph{convex} functions, that is of high interest when examining algorithm \eqref{algo1} because it seems to cover the basic idea behind \eqref{algo1}. The algorithm given in \cite{deuflhard2004newton} is:
\begin{equation}
x^{(it)}=x^{(it-1)}-s^{(it-1)}f^{(0)}(x^{(it-1)})f^{(1)}(x^{(it-1)}), ~s^{(it-1)}>0, \label{algo_origin}
\end{equation}
where $s^{(it-1)}$ denotes the step size. In \cite{deuflhard2004newton}, there is also a Lemma guaranteeing the existence of an appropriate $s^{(it-1)}$ for each iteration. There is also a discussion about a step size strategy and the convergence rate of the iterative algorithm \eqref{algo_origin} is shown to be ``linear even for bad initial guesses - however, possibly arbitrarily slow''\cite{deuflhard2004newton}, i.e., with appropriately chosen $s^{(it-1)}$ in each iteration the algorithm converges but not faster than linear. Additionally, $s^{(it-1)}$ can be chosen well enough to converge to a root but still such that the convergence of the algorithm is arbitrarily slow. In \cite{deuflhard2004newton}, it is also said that the so-called ``pseudo-convergence'' characterized by ``small'' $\Vert f^{(1)}(x)f(x)\Vert$ may occur far from the solution point due to local ill-conditioning of the Jacobian matrix. 

It is relatively easy to see that algorithm \eqref{algo1} is a special case of \eqref{algo_origin} and why it does work for convex functions (or for concave functions by changing the sign of $\alpha$).
Recall that we optimized for each possible algorithm that uses up to second derivatives. Algorithm \eqref{algo1} was discovered for the equation $x\exp(x)-1=0$ over $[-2,2]$ with starting point $x^{(0)}=0.1$. In this special case, step size $\frac{1}{f^{(2)}(x^{(it-1)})}$ was good enough to converge within the 10 allowed iterations and $\alpha$ can be set to $-1$. We could rethink the costs of a changing $\alpha$ and the usage of the inverse of the second derivative in order to force $1\neq\alpha\neq -1$.

It is interesting to consider the convergence of such algorithms for the case that the root of the function $f$ exists only in the limit. 
Then, the algorithms searching for the roots can only reach $\epsilon$-convergence.
For instance, consider $f(x)=x\exp(x)$, which is convex over, e.g., $X=(-\infty,-2]$. It holds that $\lim\limits_{x \to-\infty}f(x) = 0$ and for, e.g., the starting point $x^{(0)}=-5$, algorithm \eqref{algo1} moves in the negative $x$-direction, i.e., the algorithm iterates to the left in each iteration. 
Let us now consider the algorithm given by \eqref{algo2}. Here $\frac{f^{(1)}(x^{(it-1)})}{(f^{(2)}(x^{(it-1)}))^2}$ operates as the step size $s^{(it-1)}$. Algorithm \eqref{algo2} is more problematic, but in the special case considered here, it converges. 

All in all, $\frac{1}{f^{(2)}(x^{(it-1)})}$ and $\frac{f^{(1)}(x^{(it-1)})}{(f^{(2)}(x^{(it-1)}))^2}$ can act as step size estimators in algorithm \eqref{algo_origin} for specific functions.
Nevertheless, the value of the algorithms furnished is questionable. We show an example where $\frac{1}{f^{(2)}(x^{(it-1)})}$ and $\frac{f^{(1)}(x^{(it-1)})}{(f^{(2)}(x^{(it-1)}))^2}$ cannot act as step sizes and the algorithm would not be feasible for $\alpha=-1$. Consider the function $f(x)=x^2-3$ with $f^{(1)}(x)=2x$ and $f^{(2)}(x)=2$. Algorithms \eqref{algo1} and \eqref{algo2} both diverge for all $x\in\mathbb{R}$ except for the roots and the minimum, and some special cases where the step size $\frac{1}{f^{(2)}(x^{(it-1)})}=\frac{1}{2}$ or $\frac{f^{(1)}(x^{(it-1)})}{(f^{(2)}(x^{(it-1)}))^2}=\frac{x}{2}$ is exact, e.g., $x^{(0)}\in\{-2,2\}$. For example, considering \eqref{algo1} with the starting points $x_1^{(0)}=3$ and $x_2^{(0)}=0.1$, we get the following iteration sequence:
\begin{alignat*}{2}
  &x_1^{(0)}=3 ,~&&x_2^{(0)}=0.1 \\
  &x_1^{(1)}=-15 ,~&&x_2^{(1)}=0.399 \\
	&x_1^{(2)}=3315 ,~&&x_2^{(2)}=1.532478801 \\
	&x_1^{(3)}=-36429267622 ,~&&x_2^{(3)}=2.53090211 \\
	&x_1^{(4)}=\dots ,~&&x_2^{(4)}=\dots 
\end{alignat*}
showing divergence of algorithm \eqref{algo1} for this simple convex function $f(x)=x^2-3$.    
Algorithm \eqref{algo2} shows similar behavior. Although $f$ is convex, the algorithms both diverge showing that the step size estimators in \eqref{algo1} and \eqref{algo2} given by $\frac{1}{f^{(2)}(x^{(it-1)})}$ and $\frac{f^{(1)}(x^{(it-1)})}{(f^{(2)}(x^{(it-1)}))^2}$ cannot be optimal for all convex functions. 

\subsection{Importance of Step Size}
Recall that the optimization formulation allows choice of the step size aiming to mimic typical step size search.
In some cases this leads to ``spurious'' values of the step size $\alpha$. We allow step sizes $\alpha_{it}=\pm 2^{-\bar{\alpha}_{it}}$ with $\bar{\alpha}_{it}\in \{0,1,\dots,10\}$.
The algorithms discovered by our formulation might, in some cases, be rationalized as algorithms for optimal line search.
There are many rules for the choice of the step size when using line-search. One of those rules can be adjusted to only use step sizes of size $\alpha_{it}=\pm 2^{-\bar{\alpha}_{it}}$. The rule says then to divide the step size by 2 as long as the current step is not feasible, i.e., the step obtained in the end does not have to be optimal but only feasible.
Herein, the optimizer finds the \textit{optimal} step size for the line-search using this rule for the given function $f$ and the given iterative algorithm.
In other words, in the 1D case, step size allows each algorithm to be feasible.
Good algorithms for the general case, e.g., Newton, will have a sensible step size selection, whereas other algorithms may have an apparently spurious one.
In particular for $x \exp(x)-1=0$ the simplest algorithm described by $x^{(it-1)}+\alpha_{it}$ is the second cheapest one and the algorithm $x^{(it-1)}+\alpha_{it} f\left(x^{(it-1)}\right)^2$ is surprisingly the cheapest algorithm. 
%This is the case because the values of $f\left(x^{(it-1)}\right)^2$ have a positive effect on the direction and step size for this particular problem. 

In contrast, for two or more dimensions, such spurious behavior is not expected or obtained. A simple optimization of the step size is not enough for an algorithm to converge to a root or a minimum of a given function, since the direction provided by the gradient is crucial.
Perhaps not surprisingly, with the simple formulation not many new algorithms are found.
For instance, for the well-known and quite challenging Rosenbrock function, allowing at most 20 iterations, only Newton's algorithm is found to be feasible.
This is not surprising given the quite restrictive setting, e.g., the choice of restriction  $\alpha_{it}=\pm 2^{-\bar{\alpha}_{it}}$.
Note also that the global optimizer did not converge in all cases so that we cannot exclude the existence of feasible algorithms.
%
%\begin{itemize}
%\item Most probably, there simply are no algorithms of the form $x^{(it-1)}+\alpha_{it} \prod \left(f^{(j)}\left(x^{(it-1)} \right) \right)^{\mu_j}$ that work and are not discovered yet when allowing at most second derivatives. 
%\item The optimization instances for several algorithms may be harder than for others, e.g., when using fractions, and the optimizer needs more time to find a global solution.
%\item The algorithms need more iterations.
%\item The restriction $\alpha=\pm 2^{-\bar{\alpha}_{it}}$ is too prohibitive. 
%\end{itemize}    

\subsection{Two-Step Methods}
\label{sect:Nesterov}
We first considered the 1-dimensional strongly convex problem $\min_x \exp(x)+x^2$ over $X=[-2,2]$ with starting point $x^{(0)}=y^{(0)}=0.1$.

We also considered the minimization of $f(x_1,x_2)=(x_1-1)^2+2x_2^2-x_1x_2$ over $[-1,2]^2$ with starting point $\mathbf{x}^{(0)}=\mathbf{y}^{(0)}=(-1,-1)$. The function is convex and one could argue that it is too simple but if we want to find any algorithm that converges within a predefined maximum number of iterations $it_{con}$, there of course has to exist such an algorithm in first place. %This was not the case for our previous 2-dimensional examples. 
\begin{figure}[t]%
\centering
\begin{subfigure}{0.5\textwidth}
\begin{overpic}[width=\textwidth]{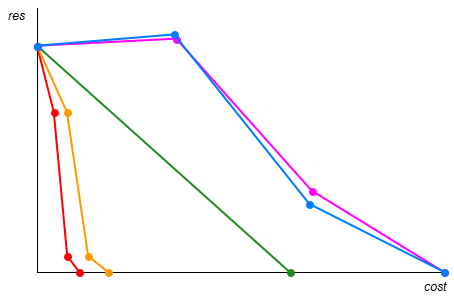}%
\put (52,40) {\begin{tikzpicture}
\draw[line width=1pt,red] (0,0)--(0.25,0);
\node[draw=none] at (1.2,0)  {\tiny$y^{(it-1)}+\alpha_{it}\nabla f$};
\draw[line width=1pt,yelloworange] (0,-0.4)--(0.25,-0.4);
\node[draw=none] at (1.4,-0.4)  {\tiny$y^{(it-1)}+\alpha_{it} f^2 \nabla f^2$};
\draw[line width=1pt,forestgreen] (0,-0.8)--(0.25,-0.8);
\node[draw=none] at (1.3,-0.8)  {\tiny$y^{(it-1)}+\alpha_{it}\frac{\nabla f}{\nabla^2 f}$};
\draw[line width=1pt,magenta] (0,-1.2)--(0.25,-1.2);
\node[draw=none] at (1.55,-1.2)  {\tiny$y^{(it-1)}+\alpha_{it}\nabla f \nabla^2 f$};
\draw[line width=1pt,cerulean] (0,-1.6)--(0.25,-1.6);
\node[draw=none] at (1.60,-1.6)  {\tiny$y^{(it-1)}+\alpha_{it}\frac{\nabla f \nabla^2 f}{f}$};
\end{tikzpicture}}
\end{overpic}
\end{subfigure}%
\begin{subfigure}{0.5\textwidth}
\includegraphics[width=\textwidth]{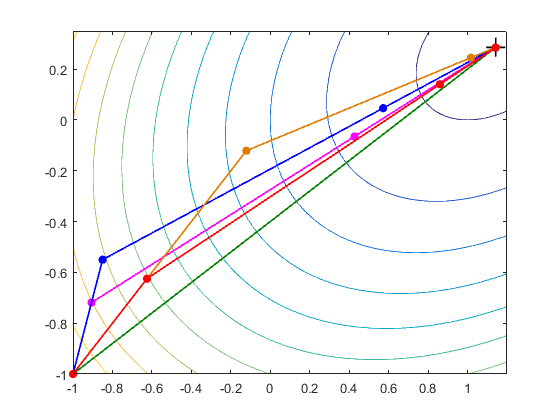}%
\end{subfigure}%
\caption{5 cheapest algorithms for the minimization of $f(x_1,x_2)=(x_1-1)^2+2x_2^2-x_1x_2$ over $X=[-1,2]^2$ with starting point $\mathbf{x}^{(0)}=\mathbf{y}^{(0)}=(-1,-1)$. Nesterov's algorithm (shown in red) is found to be the cheapest. All algorithms seem very similar as they all use gradient information with only small differences in the additional function value multiplicator and the use of hessian information.}%
\label{fig5}%
\end{figure}
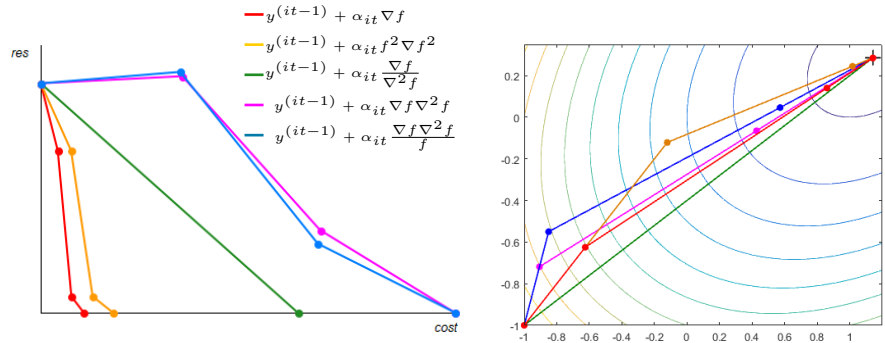
The results can be seen in Fig. \ref{fig5}. We found 11 feasible algorithms for this problem. We allowed at most 20 iterations. The 5 cheapest algorithms did not even require 5 iterations and the 6 other feasible algorithms also converged in under 5 iterations. Nesterov's method is the cheapest one found for this particular problem. In general, one could say that, except for Newton's method (green curve in Fig.~\ref{fig5}), all algorithms are alterations of Nesterov's method. 

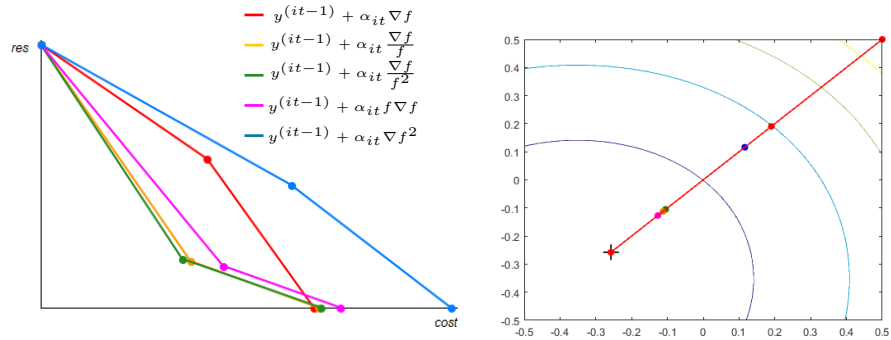
\begin{figure}[t]%
\centering
\begin{subfigure}{0.5\textwidth}
\begin{overpic}[width=\textwidth]{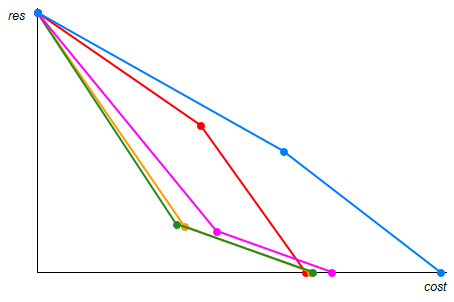}%
\put (52,40) {\begin{tikzpicture}
\draw[line width=1pt,red] (0,0)--(0.25,0);
\node[draw=none] at (1.3,0)  {\tiny$y^{(it-1)}+\alpha_{it}\nabla f$};
\draw[line width=1pt,yelloworange] (0,-0.4)--(0.25,-0.4);
\node[draw=none] at (1.3,-0.4)  {\tiny$y^{(it-1)}+\alpha_{it} \frac{\nabla f}{f}$};
\draw[line width=1pt,forestgreen] (0,-0.8)--(0.25,-0.8);
\node[draw=none] at (1.3,-0.8)  {\tiny$y^{(it-1)}+\alpha_{it}\frac{\nabla f}{ f^2}$};
\draw[line width=1pt,magenta] (0,-1.2)--(0.25,-1.2);
\node[draw=none] at (1.3,-1.2)  {\tiny$y^{(it-1)}+\alpha_{it} f \nabla f$};
\draw[line width=1pt,cerulean] (0,-1.6)--(0.25,-1.6);
\node[draw=none] at (1.3,-1.6)  {\tiny$y^{(it-1)}+\alpha_{it}\nabla f^2$};
\end{tikzpicture}}
\end{overpic}
\end{subfigure}%
\begin{subfigure}{0.5\textwidth}
\includegraphics[width=\textwidth]{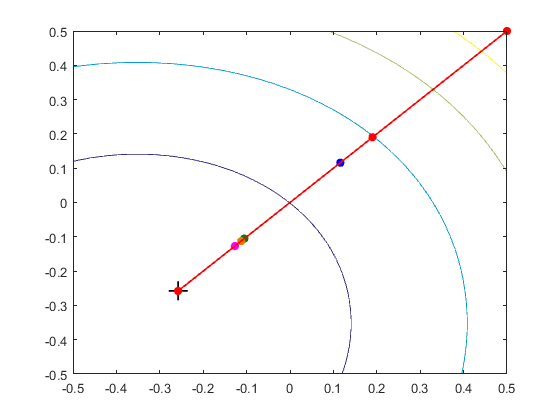}%
\end{subfigure}%
\caption{5 cheapest algorithms for the minimization of $f(x_1,x_2)=\exp(x_1)+x_1^2+\exp(x_2)+x_1^2+x_1 x_2$ over $[-1,1]^2$ with starting point $\mathbf{x}^{(0)}=\mathbf{y}^{(0)}=(0.5,0.5)$. Nesterov's algorithm is the cheapest shown in red. All algorithms seem again very similar as they all use gradient information with only small differences in the additional function value multiplicator and the use of hessian information. Note that Newton's algorithm is not among the 5 cheapest.}%
\label{fig6}%
\end{figure}

Additionally, we considered the minimization of the convex function $f(x_1,x_2)=\exp(x_1)+x_1^2+\exp(x_2)+x_1^2+x_1 x_2$ over $[-1,1]^2$ with starting point $\mathbf{x}^{(0)}=\mathbf{y}^{(0)}=(0.5,0.5)$. The optimimum of $f$ is at $\mathbf{x}^*\approx (-0.257627653,-0.257627653)$. The results can be seen in Fig.~\ref{fig6}. Note that this time Newton's algorithm is \textbf{not} among the 5 cheapest algorithms for this problem.   
It is also important to remark that the 5 algorithms all follow the same path, as can be seen in the contour plot shown in Fig.~\ref{fig6}. This is caused by the choice of the starting point and by the symmetry of the considered function.
The examples show that it is possible for us to find and discover non-trivial algorithms if the formulation is adjusted appropriately, albeit at high cost of solving the optimization formulations.
%Currently the limiting factor is obviously the computational time and the numerical bugs and errors caused by the solvers. This could become a lesser problem in near future when the research on global optimization methods and tools progresses. 

\subsection{Initial conditions}

After finding two-step algorithms for the functions mentioned in section \ref{sect:Nesterov}, we investigated the behavior of the 5 cheapest algorithms for several starting points. 

Let us first discuss the results for the strongly convex function $\exp(x)+x^2$ over $[-2,2]$ and the 5 cheapest algorithms found. We chose 17 different starting points
\[ x^{(0)}=y^{(0)} \in \{-2,-1.75,\dots,2 \}.\]
All five algorithms found for the starting point $x^{(0)}=y^{(0)}=0.1$ were also feasible for the 17 starting points considered. Nesterov's method is ``best'' among these in the sense that it was the cheapest in 15 cases and second cheapest in the remaining 2 cases. A similar statement cannot be made on the other 4 algorithms since their ranking w.r.t. cost varied depending on the initial point. Still, none of the 4 other algorithms had significant increases in cost which can probably be explained by the strong convexity of the considered function and the similarity with Nesterov's method manifested in the 4 algorithms.  

Next, let us discuss the results for $f(x_1,x_2)=(x_1-1)^2+2x_2^2-x_1x_2$ over $[-1,2]^2$ and the 5 cheapest algorithms shown in Fig.\ref{fig5} for starting point $\mathbf{x}^{(0)}=\mathbf{y}^{(0)}=(-1,-1)$. We chose 15 different starting points $\mathbf{x}^{(0)}=\mathbf{y}^{(0)} \in \{(2,2)$, (2,1), (2,0), (2,-1), (1,2), (1,1), (1,0), (1,-1), (0,2), (0,1), (0,0), (0,-1), (-1,2), (-1,1), (-1,0), $(-1,-1)\}$. Nesterov's method and Newton converged for all starting points and Nesterov was the cheapest for every initial condition. The other algorithms did not converge for all starting points or even became infeasible, e.g., the algorithm $y^{(it-1)}+\alpha_{it} f^2 \nabla f^2$ was infeasible for every starting point containing a 0 coordinate. Still, even though the unknown algorithms were infeasible for some starting points, the ranking w.r.t. the cost did not change for the cases where all algorithms converged. 

Last, let us discuss the results for $f(x_1,x_2)=\exp(x_1)+x_1^2+\exp(x_2)+x_1^2+x_1 x_2$ over $[-1,1]^2$ and the 5 cheapest algorithms shown in Fig.\ref{fig6} for starting point $\mathbf{x}^{(0)}=\mathbf{y}^{(0)}=(0.5,0.5)$. We chose 13 different starting points $\mathbf{x}^{(0)}=\mathbf{y}^{(0)} \in \{(-1,-1)$, (-1,0), (-1,1), (0,-1), (0,0), (0,1), (1,-1), (1,0), (1,1), (0.4,0.3), (0.9,-0.1), (-0.6,-0.8), $(-0.3,-0.9)\}$. Here we could observe a behavior that was not quite expected. The algorithms only converged for the each of 7 starting points $\{(-1,-1),$ (0,0), (0,1), (1,-1), (1,0), (1,1), $(0.4,0.3)\}$. Some of these 7 starting points are placed on the diagonal path taken by the algorithms, seen in Fig.\ref{fig6} and the other starting points seemed to simply provide useful derivative information. The algorithms either were infeasible for the other 6 starting points or the optimizer did not converge in the given time. Likely the infeasibility is since we used (as commonly done) the same $\alpha$ and $\beta$ for both coordinates in the algorithms and since the chosen function has symmetrical properties, the algorithms only converge under some specific conditions. In that sense the proposed method gives insight into devising new algorithmic ideas, e.g., different $\beta$ for different coordinates.

\section{Discussion of Possible Extensions}
\label{sec:extendprob}
In principle, the proposed formulation can be applied to any problem and algorithm.
In this section we list illustrative extensions to other problems of interest, starting with relatively straightforward steps and progressing towards exploratory possibilities.

\emph{Optimal Tuning of Algorithms.}
Many algorithms have tuning parameters.
Optimizing these for a given fixed algorithm using similar formulations as presented is straightforward.
%do we touch on the issue of definition of different algorithm vs algorithm with some tuning
Of course, alternative proposals exist, e.g.,~\cite{weinan2016}.

\emph{Matrix.}
Regarding the possibility of working with matrices, a formulation for finding algorithms from the family of Quasi-Newton methods could be formulated. The general idea of the Quasi-Newton methods is to update an approximate Hessian matrix with the use of gradient difference $\nabla f(x^{(it)})-\nabla f(x^{(it-1)})$. Then an approximation of the Hessian matrix is computed by the use of, e.g., Broyden's method or the Davidon-Fletcher-Powell formula, which both can be expressed with one equation.
%%This would yield a formulation of the form similar to
%%\begin{align}
%%\min \sum_{it=1}^{it_{con}}&\phi(g(x^{(it)},z^{(it)})) \label{quasi-newton}\\
%%\text{s.t. } x^{(it)}&=x^{(it-1)}+g(x^{(it-1)}), \quad i=1,\dots, it_{con} \notag\\
%%z^{(it)} &=\nabla f(x^{(it)})-\nabla f(x^{(it-1)}), \quad i=1,\dots, it_{con} \notag\\
%%H^{(it)}(f(x^{(it)})) &= \text{formula depending on } (x^{(it)},z^{(it)}),  \quad i=1,\dots, it_{con}.\notag
%%\end{align}
%%Note that $H$ is usually independent of $g$ but rather depends on $f$ and several other factors (at least for the simplest formulas) but in general a dependency of $H$ on $g$ cannot be excluded. Again this could describe potential future work as, at least to Najman, this does not seem easy at all.

\emph{Rational Polynomial Algorithms.}
A rather obvious generalization is to consider not just a single monomial but rather rational polynomials involving the derivatives.
More generally this could be extended to include non-integer and possibly even irrational powers.
This would also allow algorithms involving Taylor-series expansion.
No conceptual difficulty exists for an optimization formulation similar to the above but the number of variables increases.

\emph{Integral Operations.} The formulation essentially also allows algorithms that perform integration, if $f^{(j)}$ with $j<0$ is allowed.
Obviously for multivariate programs ($n_x>1$) the dimensionality needs to be appropriately considered.

\emph{Implicit Methods.} Implicit methods, e.g., implicit Euler, do not only involve evaluation of derivatives but also the solution of a system of equations at each step. %, e.g., in the implicit Euler
%\[ x^{(it)}=x^{(it-1)}+f(x^{(it)}) \Delta t.\]
In other words we cannot directly write such methods as monomials in Definition~\ref{def:ratpolalgo}.
However, it is conceptually relatively easy to adjust the formulation, in particular by changing the update scheme to something along the lines of
$$  x^{(it)}=x^{(it-1)} + \alpha \Pi_{j=0}^{j_{max}} \left(f^{(j)}(x^{(it)})\right)^{\nu_{j,1}}. $$
Computationally, the optimization problems will become substantially harder.
The variables $x^{(it)}$ cannot be directly eliminated but rather the equations have to be given to the optimizer or addressed as implicit functions, see~\cite{barton_06_1,stuber_14_1}.

\emph{General Multistep Methods.} Explicit multistep methods are those that calculate the current iterate based not only on the immediately previous but rather also on additional preceding steps.
Examples include multistep integration methods as well as derivative-free equation solution methods such as secant or bisection.
It is straightforward to extend the formulations to allow such methods.
e.g., for two-step methods
\begin{align*}
  x^{(it)}=x^{(it-1)} &+ \alpha_1 \Pi_{j=0}^{j_{max,1}} \left(f^{(j)}(x^{(it-1)})\right)^{\nu_{j,1}}\\
  &+ \alpha_2 \Pi_{j=0}^{j_{max,2}} \left(f^{(j)}(x^{(it-2)})\right)^{\nu_{j,2}}.
  \end{align*}
Obviously the number of optimization variables increases accordingly.
This can also be used for methods that rely on approximation of derivatives, e.g., Quasi-Newton methods.
Similarly, methods such as secant can be encoded by allowing mixed polynomials of ${\bf x}$ along with ${\bf f}^{(j)}$.
Bisection is not captured as this requires if-then-else; however, it should be relatively easy to capture this using integer variables which fits the proposed MINLP framework.
Obviously for \emph{implicit multistep} methods the two preceding formulations can be combined
with a formulation along the lines of
\[ x^{(it)}=x^{(it-1)}+ \sum_{\hat{it}=\hat{it}_{min}}^{\hat{it}_{max}} \alpha_{it,\hat{it}} \Pi_{j=0}^{j_{max,1}} \left(f^{(j)}(x^{(it+\hat{it})})\right)^{\nu_{j,\hat{it}}}. \]

\emph{Global Optimization.}
A challenge in global optimization is that there are no explicit optimality criteria that can be used in the formulation.
The definition of global optimality $ f\left({\bf x}^*\right) \leq f\left({\bf x}\right)$, $\forall {\bf x} \in X $
can be added to the optimization of algorithm formulation but this results in a (generalized) semi-infinite problem (SIP).
There are methods to deal with SIP problems but they are computationally very expensive~\cite{mitsos_08_siprestricted}.
Another challenge is that deterministic and stochastic global methods do not follow the simple update using just monomials as in Definition~\ref{def:ratpolalgo}.
As such, major conceptual difficulties are expected for such algorithms, including an estimation of the cost of calculating the relaxations.
Moreover, to estimate the cost of such algorithms the results of ongoing analysis in the spirit of automatic differentiation will be required~\cite{du_94_1,mitsos_10_mcrate}.
and the cost does not only depend on the values of the derivatives at a point.
For instance in $\alpha$BB~\cite{maranas_92_1,adjiman_96_1} one needs to estimate the eigenvalues of the Hessian for the domain; in McCormick relaxations~\cite{mccormick_76_1,mitsos_12_mcmul} the cost depends on how many times the composition theorem has to be applied.
Recall also the discussion on cost for some local algorithms such as line-search methods.

\emph{Optimal Algorithms with Optimal Steps.}
Another interesting class are algorithms that include an optimization step within them, such as selection of an optimal step, e.g., \cite{zhang2014spectral,wibisono2016variational}. % ,wibisono2015,lee2016gradient
%Each of the algorithms depends on calculating an optimal step/gradient in the sense of finding a $\min/\max$ in order to construct the next iteration step.
%Such an algorithm family consists of the generalized Gauss-Newton methods \cite{bock1983recent}.
We can consider a bilevel formulation of the form
\begin{align}
\min \sum_{it=1}^{it_{con}}&\phi(g(x^{(it)},z^{(it)})) \label{bilevel}\\
\text{s.t. } x^{(it)}&=g(x^{(it-1)},z^{(it-1)}), \quad i=1,\dots, it_{con} \notag\\
z^{(it)} &\in \arg \min_{z'} f(x^{(it-1)},z'), \quad i=1,\dots, it_{con}, \notag
\end{align}
describing problems where the next step is given by an optimality condition, e.g., $f$ could describe the least squares formulation for the substep in a generalized Gauss-Newton algorithm in which the outer problem is similar to the formulation above, while the lower-level problem uses a different objective to calculate the optimal step.
There exist (quite expensive) approaches to solve bilevel programs with nonconvexities~\cite{mitsos_06_2}.
When the lower level problem is unconstrained regarding $f$ can be avoided, as is the case in generalized Gauss-Newton algorithms, then the problem becomes much simpler.
Note that although we have several $\arg \min$, they are not correlated, i.e., we do not obtain a $n$-level optimization problem but rather a bilevel problem with multiple lower level problems.
With such a formulation we might be able to prove optimality of one of the three gradient methods described in \cite{wibisono2016variational} or discover alternatives.
Additionally it would be possible to discover well-known methods such as the nonlinear conjugate gradient method \cite{luenberger1973introduction}, where the line search is described by the $\arg \min$ operator or even stochastic methods which often depend on the computation of the maximum expected value of some iterative random variable $\mathcal{X}^{(it)}$.
Similarly, it could be possible
%to accordingly adjust the above formulation \eqref{bilevel}
to find the desired controller algorithm in \cite{economou1985}.

\emph{Continuous Form.}
%{\tt see Yannis emails 05.04.16 in exchange with CW Gear}
Discrete methods are often reformulated in a continuous form.
For instance, Boggs proposed a continuous form of Newton's method~\cite{boggs_71_1,boggs_76_1}
$ \dot{x}(t)=\frac{f(x(t))}{f^{(1)}(x(t))}$.
See also the recent embedding of discrete algorithms like Nesterov's scheme in continuous implementations \cite{jordan_16_1,candes_15_1}.
It seems possible to consider the optimization of these continuous variants of the algorithms using similar formulations.
%To capture the monomial algorithms, Definition~\ref{def:ratpolalgo}, integer variables can be used.
%Allowing for a change of the algorithm (similar to $\nu_{j}^{it_1} \neq \nu_{j}^{it_2}$ in the discrete case), these integer variables need to be dependent on $t$.
Typically discretization methods are used to optimize with such dynamics embedded.
%Moreover, mixed-integer dynamic optimization problems are notoriously hard to solve.
Thus, this continuous formulation seems more challenging to solve than the discrete form above. % in~\ref{sec:optformufin}
If a particular structure of the problem can be recognized, it could however be interesting, for instance to apply a relaxation similar to~\cite{sager_09_1}.

\emph{Dynamic Simulation Algorithms.}
The task here is to simulate a dynamical system, e.g., ordinary differential equations (ODE)
$ \dot{\bf x}(t)={\bf f}({\bf x}(t)) $ along with some initial conditions ${\bf x}(t=0)={\bf x}^{init}$ for a given time interval $[0,t_f]$.
We need to define what is meant by a feasible algorithm.
A natural definition involves the difference of the computed time trajectory from the exact solution, which is not known and thus does not yield an explicit condition.
One should therefore check the degree to which the resulting approximate curve, possibly after interpolation, satisfies the differential equation over this time interval.

\emph{Dynamic Optimization Algorithms.}
Dynamic optimization combines the difficulties of optimization and dynamic simulation.
Moreover, it results in a somewhat amusing cyclic problem: we require an algorithm for the solution of dynamic optimization problems to select a dynamic optimization algorithm.
As aforementioned, this is not prohibitive, e.g., in the offline design of an algorithm to be used online.

\emph{Algorithms in the Limit.}
Considering algorithms that provably converge to the correct solution in the limit makes the optimization problem more challenging.
The infinite-dimension formulation \eqref{optformuinf} is in principle applicable with the aforementioned challenges.
If the finite (parametrized) formulation \eqref{optformufin}, was directly applied, an infinite number of variables would have to be solved for.
In such cases one could test for asymptotic self-similarity of the algorithm behavior as a way of assessing its asymptotic result.

\emph{Quantum Algorithms.}
A potential breakthrough in computational engineering would be realized by quantum computers.
These will require new algorithms, and there are several scientists that are developing such algorithms.
It would thus be of extreme interest to consider the extension of the proposed formulation to quantum algorithms and/or their real-world realizations.
This may result in ``regular'' optimization problems or problems that need to be solved with quantum algorithms themselves.

\section{Limitations}
The proposed formulation is prototypical and a number of challenges arise naturally.
As aforementioned, the knowledge of an explicit cost may not be a good approximation for all problems.
It is also clear that different objective functions, including, e.g., considerations of memory usage or the inclusion of error correction features, 
may dramatically effect the results.

We expect the proposed optimization formulation to be very hard to solve and the numerical results confirm this.
In particular we expect it to be at least as hard as the original problems.
Proving that statement is outside the scope of the manuscript but two arguments are given, namely that the final constraint corresponds to solution of problem and that the subproblems of the optimization problem are the same or at least related to the original problem.

Herein brute force and general-purpose solvers are used for the solution of the optimization problems.
It is conceivable to develop specialized numerical algorithms that will perform much better than general-purpose ones due to the specific structure.
In essence we have an integer problem with linear objective and a single nonlinear constraint.
%A challenge could be that for noninteger values of the variables (as is done in normal constraints) the operations do not seem well defined.
It seems promising to mask the intermediate variables from the optimizer. This would in essence follow the optimization with implicit functions embedded \cite{barton_06_1} and follow-up work.
It is also conceivable to move to parallel computing, but suitable algorithms are not yet parallelized.

In our formulation, $f$ is assumed to be a given function, so that we can apply current numerical methods for the optimization.
It is, however, possible to consider multiple instances of the problem simultaneously, i.e., allow $f$ to be among a class of functions.
This can be done similar to stochastic optimization~\cite{birge_book_1}.
A simple method is to sample the instances of interest (functions $f$) and optimize for some weighted/average performance of the algorithm.
Alternatively, the instances can be parametrized by continuous parameters and the objective in the above formulations replaced with some metric of the parameter distribution, e.g., the expectation of the cost.
It is also possible, and likely promising, to consider worst-case performance, e.g., in a min-max formulation~\cite{falk_77_1}. 
It is also interesting to consider the relation of this work to the theorems developed in \cite{wolpert1997no},
in particular that optimality of an algorithm over one class of problems does not guarantee optimality over another class.

\section{Relation to Computational Complexity}
As mentioned above, we expect our formulation to be at least as hard as the original problems.
This would imply that the formulation cannot be easy to solve when the problems themselves are NP-hard.

A potential breakthrough would be if for a particular instance of a given NP-hard problem an algorithm is furnished that is shown to have polynomial complexity. This would directly imply that P=NP.
In contrast, if algorithms of exponential complexity are furnished for a given instance problem this does not prove that P $\neq$ NP, although the conclusion is tempting (``if the best-possible algorithm is exponential, then no polynomial algorithm exists'').
To begin with, all possible algorithms would have to be allowed in the formulation to be able to conclude that no better algorithm exists.
%
%  the issue of open/closed sets of searches comes back to me- is the limit of what we find one of the things we find ----
%
Moreover, we propose methods to solve for a \emph{fixed size} problem.
The best algorithm for that size may exhibit exponential complexity, but this does not prove that there does not exist a polynomial algorithm which will be better for very large sizes. 
For instance, consider linear programs: simplex variants with exponential complexity outperform interior-point methods of polynomial complexity for all but the largest sizes.
It is however, conceivable to use our formulation to obtain some insight. For instance, solving instances of different sizes and observing if the same algorithm is found optimal or if the optimal algorithm depends on the instance size.
% We should also consider the approach presented by \cite{Macready:1995:MOP:228587.228601} where the authors focus on a particular problem and analyze occurring effects over all possible algorithms. In our case we %are able to analyze effects over several possible iterative algorithms.

\section{Conclusions}
An MINLP formulation is proposed that can readily devise optimal algorithms (among a relatively large family of algorithms)
for several prototypical problems, including solution of nonlinear equations and nonlinear optimization.
Simple numerical case studies demonstrate that well-known algorithms can be identified along with new ones.
We argue that the formulation is conceptually extensible to many interesting classes of problems, including quantum algorithms.
Substantial work is now needed to develop and implement these extensions so as to numerically devise optimal algorithms for interesting classes of challenging problems where such algorithms are simply not known.
Also, the similarity to model-reference adaptive control (MRAC) and internal model control (IMC) can be further explored in the future.
The optimization problems obtained can, however, be very challenging and no claim is made that optimal algorithm discovery will be computationally cheap.
However, in addition to the theoretical interest, there are certainly applications.
Finding guaranteed optimal algorithms for a given problem implies understanding/classifying the difficulty of this problem.
And it will certainly be worthwhile to automatically devise algorithms offline for problems to be solved online.

\appendix

\section{Solution of Equation Systems}
For illustration purposes we will discuss the solution of equation systems using only zero and first derivatives.
The zero derivative $f^{(0)}\left({\bf x}^{(it)}\right)$  is a vector of the same size as the point ${\bf x}$.
By definition, the required polynomial expression of the derivative $\left(f^{(0)}\left({\bf x}^{(it)}\right)\right)^{\nu_k}$ is calculated element by element.
Thus, we can equivalently rewrite each element $if$ of the vector $\left({\bf f}^{(0)}\left({\bf x}^{(it-1)}\right)\right)^{\nu_k}$ as $\sum_{k=k_{min}}^{k_{max}} y_0^k \left(f^{(0)}_{if}\left({\bf x}^{(it-1)}\right)\right)^k$ along with a constraint $\sum_{k=k_{min}}^{k_{max}} y_0^k=1$.
In contrast, the first derivative is the Jacobian matrix ${\bf J}$ with elements $\left.\frac{\partial f_{if}}{\partial x_{ix}}\right|_{{\bf x}^{(it)}}$.
In principle the inverse can be written analytically but this is only realistic for small-scale systems.
For systems of substantial sizes, it is better to calculate the direction via the solution of an equation system.
\[  {\bf J}\left({\bf x}^{(it-1)}\right) \Delta {\bf x}= \left({\bf f}^{(0)}\left({\bf x}^{(it-1)}\right)\right)^{\nu_k},\]
the solution of which gives us the direction $\Delta {\bf x}$ to be taken if the optimizer selects the inverse of the Jacobian, i.e., sets $\nu_1=-1$.
Similarly for other values of $\nu_1$ we have to calculate other steps via explicit matrix multiplication or via similar equation systems.
As aforementioned, the equations can be explicitly or implicitly written, see~\cite{barton_06_1} and its extensions~\cite{stuber_14_1}.

\section{Alternative to Obtain a Regular Optimization Problem}
For the sake of simplicity an alternative is presented to obtaining a finite number of variables in the optimal control problem.
Again a maximal number of iterations $it_{max}$ is considered and the cost is penalized by the residual, e.g.,
\begin{align}
&\min_{g, {x}^{(it)}} \sum_{it=1}^{it_{max}} \left|f\left(x^{(it)}\right)\right| \phi\left(g\left(x^{(it)}\right)\right) \notag \\
&\text{s.t. } x^{(it)}=g\left(x^{(it-1)}\right), \quad it=1,2,\dots, it_{max}  \\
& \left|f\left(x^{(it_{max})}\right)\right| < \varepsilon,
\end{align}
where again the last constraint ensures convergence.
Without it, if algorithms with zero cost exist, then these will be chosen even if they are infeasible; they will give a zero objective value without meeting convergence.

\section{Obtaining a Regular MINLP}
The finite problem proposed contains terms $\left(f^{(j)}(x^{(it-1)})\right)^{\nu_j}$.
We thus exponentiate a potentially negative real-valued ($f^{(j)}(x^{(it-1)})$) to the power $\nu_j$.
At feasible points  $\nu_j$ is integer and thus the operation well-defined.
However, in intermediate iterations of the optimization, e.g., in branch-and-bound, $\nu_j$ may be real-valued and thus the operation is not defined.
Consequently, the optimization formulation cannot be directly solved using standard MINLP solvers.
It is however, relatively easy to reformulate as an MINLP with linear dependence on the integer variables.
One way is to introduce for each $j$ as many binary variables $y_j^k$ as we have potential values for $\nu_j$ and enforce that exactly one is zero (special case of SOS-1 set).
Then the term $f^{(j)}(x^{(it-1)})$ is equivalently rewritten as $\sum_{k=k_{min}}^{k_{max}} y_j^k \left(f^{(j)}(x^{(it-1)})\right)^k$ along with a constraint $\sum_{k=k_{min}}^{k_{max}} y_j^k=1$.

\section*{Acknowledgments}
IGK and AM are indebted to the late C.A. Floudas for bringing them together. Fruitful discussions with G.A. Kevrekidis and the help of C.W. Gear with the continuous formulation are greatly appreciated. 

This project has received funding from the Deutsche Forschungsgemeinschaft (DFG, German Research Foundation) ``Improved McCormick Relaxations for the efficient Global Optimization in the Space of Degrees of Freedom'' MA 1188/34-1. The work of IGK was partially supported by the US National Science Foundation (CDS\&E program), by the US AFOSR (Dr. F. Darema) and DARPA contract  HR0011-16-C-0016.

% Bibliography
\bibliographystyle{abbrv}
{\footnotesize
\bibliography{optimalalgorithms}}

\end{document}